\documentclass[reprint,superscriptaddress,nofootinbib,amsmath,amssymb]{revtex4-1}
\usepackage{xcolor}
\usepackage{graphicx} 
\usepackage{float}
\newtheorem{thm}{THEOREM}[section]

\newtheorem{prp}[thm]{PROPOSITION}

\newcommand{\refeq}[1]{(\ref{#1})}

\newcommand{\vect}[1] {\boldsymbol{{ #1}} }

\newcommand{\Rset}{\mathbb{R}}

\newcommand{\aV}{\vect{a}}              

\newcommand{\sV}{{\vect{s}}}            

\newcommand{\NullV}{\vect{0}}

\newcommand{\abs}[1]{\big\vert #1 \big\vert}

\renewcommand{\leq}{\leqslant}
\renewcommand{\geq}{\geqslant}

\newcommand{\Efrak}{{\mathfrak E}} 
\newcommand{\Ffrak}{{\mathfrak F}} 

\begin{document}
\title{On the asymptotic decay of the Schr\"odinger--Newton ground state} \vspace{-10pt}
	
\author{Michael K.-H. Kiessling} 

\affiliation{Department of Mathematics, Rutgers University,
                110 Frelinghuysen Rd., Piscataway, NJ 08854, USA}
\email{$\phantom{.}$\vspace{-20pt}\copyright (2021) The author.\vspace{-32pt}}

\begin{abstract}
\noindent 
 The asymptotics of the ground state $u(r)$ of the Schr\"odinger--Newton equation in $\Rset^3$
was determined by V. Moroz and J. van Schaftingen to be $u(r) \sim A e^{-r}/ r^{1 - \|u\|_2^2/8\pi}$ for some $A>0$,
in units in which the ground state energy is $-1$. 
  They left open the value of $\|u\|_2^2$, the squared $L^2$ norm of $u$.
 Here it is rigorously shown that $2^{1/3}3\pi^2\leq \|u\|_2^2\leq \textcolor{red}{8}\pi^{3/2}$. 
 It is reported that numerically $\|u\|_2^2\approx 14.03\pi$, revealing
that the monomial prefactor of $e^{-r}$ increases with $r$ in a concave manner.
 Asymptotic results are proposed for the Schr\"odinger--Newton equation with external $\sim - K/r$ potential, and 
for the related Hartree equation of a bosonic atom or ion. 
\vspace{-15pt}
\end{abstract}
$\phantom{xi}$\hfill 1 
\maketitle

 $\phantom{nix}$\vspace{-1.7truecm}

\section{Introduction}  \vspace{-10pt}

 The non-linear integro-differential equation \vspace{-5pt}
\begin{equation}
-  \Delta u(\sV) + u(\sV) =
 \displaystyle{\frac{1}{4\pi}\int_{\Rset^3} \frac{1}{\abs{\sV-\sV'}}|u|^2(\sV') d^3\!s'}\, u(\sV)\vspace{-5pt}
 \label{CPSN}
\end{equation}
shows up in a variety of models in physics, and as such
is variously known as Pekar's equation, Choquard's equation, Schr\"odinger--Newton equation, and by other names as well. 
 Pekar \cite{Pek} suggested it as an approximate model for H. Fr\"ohlich's condensed matter polaron \cite{FrA}, \cite{FrB};
subsequently it was derived in a suitable limit $N\to\infty$ from $N$-body QM, see \cite{DVa}, \cite{DVb}, \cite{LT}.
Choquard proposed it to
characterize a self-trapped electron in a one-component quantum plasma; cf. \cite{LiebC}, \cite{ChStV}.
 The name Schr\"odinger--Newton equation was coined by R. Penrose in his proposal that quantum-mechanical wave function
collapse is caused by gravity \cite{Pen}. 
 It also occurs in the theory of hypothetical bosonic stars in the context of the dark matter mystery; it can
be derived from QM of $N$ gravitating spin-zero bosons by taking a Hartree limit $N\to\infty$, see \cite{Lewin}
and references therein.
 A survey of rigorous results for equation \refeq{CPSN} and some of its generalizations is in \cite{MvSreview}.
   \vspace{-10pt}

 Positive solutions to \refeq{CPSN} minimize the functional \vspace{-10pt}
\begin{widetext}\vspace{-10pt}
\begin{eqnarray} \label{eq:uFctl}
\Efrak(u) := 
\int_{\Rset^3} \abs{\nabla u}^2(\sV) d^3\!s  
+
\int_{\Rset^3} |u|^2(\sV) d^3\!s 
- \frac{1}{8\pi} \int_{\Rset^3} \int_{\Rset^3} 
\frac{|u|^2(\sV)|u|^2(\sV')}{\abs{\sV-\sV'}}d^3\!s d^3\!s'
\end{eqnarray}\vspace{-10pt}
\end{widetext}
$\phantom{nix}$\vspace{-25pt}

\noindent
over the set of $u\in H^1(\Rset^3)$ under the Nehari condition that $u\neq 0$ satisfies
${4\pi} \int_{\Rset^3}\big( \abs{\nabla u}^2(\sV) + |u|^2(\sV) \big) d^3\!s 
= 
\int_{\Rset^3} \int_{\Rset^3} {\abs{\sV-\sV'}}^{-1}{|u|^2(\sV)|u|^2(\sV')}d^3\!s d^3\!s'$
(obtained by multiplying \refeq{CPSN} with $u$ and integrating over $\Rset^3$); this is
equivalent to the usual minimization of r.h.s.\refeq{eq:uFctl} without $\int |u|^2 (\sV)d^3s$,
yet under a normalization condition on this integral.
 It is known \cite{LiebC}, \cite{MvSreview} that any positive minimizer is radially symmetric about and decreasing away from an 
arbitrary point $\sV_0\in\Rset^3$, and that for each $\sV_0$ there is a unique such solution.
 Moreover, since \refeq{CPSN} and \refeq{eq:uFctl} are invariant under translations $\sV\mapsto \sV + \aV$, without loss of 
generality we assume that $\sV_0 = \NullV$, and
(abusing notation) write $u(r)$ for $u(\sV)$ to denote this radially symmetric positive solution.
 We note that $u^\prime(0)=0$.
 
 Interestingly, the somewhat delicate question of the asymptotic decay of positive minimizers has not yet been sorted 
out completely, as far as we can tell.
 D. Kumar and V. Soni in \cite{KS} claimed that there exists a positive constant $A$ such that (in our units)
$u(r) = A e^{-r} +$l.o.t., where ``l.o.t.'' means ``lower-order terms.'' 
 Based on their asymptotic analysis they concluded that (in units in which $\|u\|_2^{}=1$) the energy coincides
with that of the hydrogen atom, but a nonlocal quantity like an eigenvalue cannot be determined by a truncated
asymptotic expansion, and indeeed their energy claim was subsequently disproved in \cite{T}.
  K. P. Tod and I. Moroz in \cite{TM} in turn claimed, though without proof, that if $u(r)$ is a positive radial solution 
to \refeq{CPSN}, then $u(r) = A e^{-r}/r +$l.o.t.; cf (2.17a/b) in \cite{TM}. 
 This claim was announced in \cite{MPT} and repeated in \cite{Hthesis}.
 What \emph{is} proved in \cite{TM} is that for every positive $C<1$ 
there are $A>0$ and $b>0$ such that $u(r) < Ae^{-Cr}/r$ for all $r>b$, see Thm.3.1 in \cite{TM}
    (incidentally, a factor $e^{Cb}$ is obviously missing at r.h.s.(3.9) of \cite{TM}),
but such an upper bound alone cannot establish the asymptotic behavior claimed in \cite{TM}.
 The question of the asymptotic behavior of $u(r)$ was taken up again by V. Moroz and J. van Schaftingen \cite{MvS},
who proved that the unique radially symmetric positive solution $u(r)$ to \refeq{CPSN} obeys 
$\lim_{r\to\infty} r^{1 - \|u\|_2^2/8\pi} e^r u(r) \in (0,\infty)$; thus one has the asymptotic law
$u(r) = A e^{-r}/ r^{1 - \|u\|_2^2/8\pi}+$l.o.t. for some $A>0$; see their Thm.4, case $p=2$ with $N=3$ and $\alpha=2$;
and see also section 3.3.4 in \cite{MvSreview}.
 However, the value of $\|u\|_2^{2}$ was left open in \cite{MvS}.
 While the results of \cite{MvS} rule out the asymptotic behavior claimed in \cite{TM}, \cite{MPT}, and \cite{Hthesis},
it does leave open the possibility that the ground state solution of \refeq{SN} perhaps decays as claimed in \cite{KS},
i.e. purely exponentially like the hydrogenic ground state (this would be the case if $\|u\|_2^2 = 8\pi$); but
the exponential function could also have a decaying monomial factor (i.e. $\|u\|_2^2 < 8\pi$) or an increasing 
one (i.e. $\|u\|_2^2 > 8\pi$); in the latter case: is it concave 
\textcolor{blue}{(i.e. $\|u\|_2^2 < 16\pi$)}, linear \textcolor{blue}{($\|u\|_2^2 = 16\pi$)}, 
or convex \textcolor{blue}{($\|u\|_2^2 > 16\pi$)}? 
\newpage

 With the help of rigorous energy bounds similar to those in \cite{LiebC} or \cite{T}, we prove
\begin{prp}\label{prop} 
The $L^2$ norm of the positive $H^1$ solution $u(r)$ of \refeq{SN} obeys the bounds
\begin{equation}\label{normOFuSQRbounds}
2^{1/3}3\pi^2\leq \|u\|_2^2\leq 
\textcolor{red}{8}\pi^{3/2}. 
\end{equation}
\end{prp}
 This is strong enough to rigorously rule out the asymptotic form of $u(r)$ proposed in \cite{KS} (since $2^{1/3}3\pi^2 > 8\pi$),
showing that the monomial prefactor of $e^{-r}$ is increasing with $r$,
\textcolor{red}{and strong enough (since $8\pi^{3/2}<16\pi$) to prove that
the monomial prefactor of $e^{-r}$ increases in a strictly concave manner.}

 Although numerical studies of the Schr\"odinger--Newton equation have been carried out \cite{MPT}, \cite{Hthesis}, \cite{GW},
 we are unaware of any which has addressed itself to the power of the radial monomial correction factor to the exponential 
function.
 Yet information about $\| u\|_2^2$ can be extracted from  numerical data in \cite{GW} by rescaling,
revealing that the monomial prefactor of $e^{-r}$ is $\propto r^{\beta}$ with $\beta\approx 0.754$, 
i.e. increasing and strictly concave. 
 We have carried out our own numerical study for \refeq{SN} and directly computed that $\| u\|_2^2 \approx 14.03\pi$, 
compatible with the result extracted from \cite{GW} by rescaling.

Our upper bound in \refeq{normOFuSQRbounds} is 
a factor      
\textcolor{red}{$\approx 1.011$} too large compared to the numerically computed value.

\section{Rigorous bounds on $\|u\|_2^2$}  \vspace{-10pt}

In the following we prove Proposition \ref{prop}.

\emph{Proof}: For the pupose of our proof, and also for later convenience, 
we rescale \refeq{CPSN} into 
\begin{eqnarray}\hspace{-5pt}
-  \Delta \psi(\sV) -  \displaystyle{ 2 \int_{\Rset^3} \frac{1}{\abs{\sV-\sV'}}|\psi|^2(\sV') d^3\!s'}\, \psi(\sV)
= E \psi(\sV) ,
 \label{CPSNnat}
 \vspace{-5pt}
\end{eqnarray}
in which form the Schr\"odinger--Newton equation appears in \cite{LiebC} and \cite{GW}; here, $\|\psi\|_2^2=1$.
 Equation \refeq{CPSNnat} is the Euler--Lagrange equation for the minimization of the functional
\begin{eqnarray} \label{eq:psiFctl}
\hspace{-5pt}\Ffrak(\psi) := \!\!
\int_{\Rset^3}\!\! \abs{\nabla \psi}^2(\sV) d^3\!s  
-
\! \int_{\Rset^3}\! \int_{\Rset^3} \!\!\!\!
\frac{|\psi|^2(\sV)|\psi|^2(\sV')}{\abs{\sV-\sV'}}d^3\!s d^3\!s'\
\end{eqnarray}
over the Sobolev space $H^1(\Rset^3)$ under the constraint $\int_{\Rset^3} |\psi|^2(\sV) d^3\!s =1$;
the eigenvalue $E$ is the Lagrange multiplier for this constraint; see  \cite{LiebC}.

 Let $\psi_1(r)$ be the minimizer, and for real $\lambda >0$ define $\psi_\lambda(r):= {\lambda}^{3/2}\psi_1(\lambda r)$.
 Then $\forall\,\lambda:\, \|\psi_\lambda\|_2 = 1$. 
 By noting that $\frac{d}{d\lambda}\Ffrak(\psi_\lambda)|^{}_{\lambda=1}=0$ we obtain the virial identity
\begin{eqnarray} \label{eq:virial}
\hspace{-5pt} 
2\int_{\Rset^3}\!\! \abs{\nabla \psi_1}^2(\sV) d^3\!s  
=
\! \int_{\Rset^3}\! \int_{\Rset^3} \!\!\!\!
\frac{|\psi_1|^2(\sV)|\psi_1|^2(\sV')}{\abs{\sV-\sV'}}d^3\!s d^3\!s'.\
\end{eqnarray}
 On the other hand, setting $\psi=\psi_1$ in \refeq{CPSNnat}, then 
multiplying \refeq{CPSNnat} by $\psi_1$ and integrating over $\Rset^3$, one obtains 
 \begin{eqnarray} \label{eq:energy}
\hspace{-10pt} 
E_0\! =\!\! \int_{\Rset^3}\!\! \abs{\nabla \psi_1}^2(\sV) d^3\!s  
-
2 \! \int_{\Rset^3}\!\! \int_{\Rset^3} \!\!\!\!
\frac{|\psi_1|^2(\sV)|\psi_1|^2(\sV')}{\abs{\sV-\sV'}}d^3\!s d^3\!s'\
\end{eqnarray}
for the ground state energy $E_0$.
 Now using \refeq{eq:virial} in \refeq{eq:energy} and also in \refeq{eq:psiFctl}, by comparison we obtain
 \begin{eqnarray} \label{eq:energyF}
E_0 = 3\,\Ffrak(\psi_1).
\end{eqnarray}
 Thus, any upper or lower bounds on $\Ffrak(\psi_1)$ over $H^1$ under the normalization constraint $\|\psi_1\|_2=1$
translate into corresponding upper and lower bounds on the ground state energy $E_0$.

 Next, by Sobolev's inequality (cf. \cite{T}, p.174),
 \begin{eqnarray} \label{eq:Sob}
\int_{\Rset^3}\!\! \abs{\nabla \psi}^2(\sV) d^3\!s  
\geq 
3\left[\frac{\pi}{2}\right]^{\!\frac43}
 \left( \int_{\Rset^3}|\psi|^6(\sV)d^3\!s\right)^{\!\!\frac13}.
\end{eqnarray}
 On the other hand, the Hardy--Littlewood--Sobolev inequality yields
\begin{eqnarray} \label{eq:HLS}
\hspace{-10pt}
 \int_{\Rset^3\times\Rset^3} \!\!\!\! \!\!\!\!
\frac{|\psi|^2(\sV)|\psi|^2(\sV')}{\abs{\sV-\sV'}}d^3\!s d^3\!s'
\leq
\frac43\!\left[\frac{8}{\pi}\right]^{\!\frac13}\!\!\!
 \left( \int_{\Rset^3}\!\!\!|\psi|^{\frac{12}{5}}(\sV)d^3\!s\!\!\right)^{\!\!\frac{5}{3}}\!\!\!\!,\ \ 
\end{eqnarray}
and H\"older's inequality gives (cf. \cite{T}, p.175)
 \begin{eqnarray} \label{eq:Hoe}
\hspace{-10pt}
\int_{\Rset^3}\!\! \abs{\psi}^{\frac{12}{5}}(\sV) d^3\!s  
\leq 
\left( \int_{\Rset^3}|\psi|^2(\sV)d^3\!s\!\right)^{\!\!\frac{9}{10}}\!\!\!
\left( \int_{\Rset^3}|\psi|^6(\sV)d^3\!s\!\right)^{\!\!\frac{1}{10}}\!\!\!\!\!\!,\ \ 
\end{eqnarray}
which simplifies because $ \int_{\Rset^3}|\psi|^2(\sV)d^3\!s=1$.
 Now setting $\int_{\Rset^3}|\psi_1|^6(\sV)d^3\!s=:x^{6}$, our chain of inequalities yields
 \begin{eqnarray} \label{eq:FlowerBOUNDx}
\Ffrak(\psi_1) \geq 
3\left[\frac{\pi}{2}\right]^{\!\frac43} x^2 - 
\frac43\!\left[\frac{8}{\pi}\right]^{\!\frac13} x;
\end{eqnarray}
and since we do not know $x$, to be on the safe side we minimize r.h.s.\refeq{eq:FlowerBOUNDx} w.r.t. $x>0$ and obtain
\begin{eqnarray} \label{eq:FlowerBOUND}
\Ffrak(\psi_1) \geq 
-\frac{32}{27}\frac{2^{1/3}}{\pi^2} \approx -0.1513.
\end{eqnarray}
 This yields 
\begin{eqnarray} \label{eq:ElowerBOUND}
E_0 \geq 
-\frac{32}{9}\frac{2^{1/3}}{\pi^2} \approx -0.4539.
\end{eqnarray}
 By rescaling the ground state energy to $E=-1$ we find $\|u\|_2^2 = 8\pi/\sqrt{|E_0|}$, and so, since $E_0<0$, by
\refeq{eq:ElowerBOUND} we have
\begin{eqnarray} \label{eq:L2uLowerBOUND}
\|u\|_2^2 
\geq 
 2^{1/3}3\pi^2 \approx 37.305.
\end{eqnarray}
 This proves the lower bound in Proposition \ref{prop}.

 To obtain the upper  bound in Proposition \ref{prop} we insert the Gaussian trial wave function 
$\psi_G(r) := \exp(-r^2/2R^2)/(\pi^{3/4}R^{3/2})$ into $\Ffrak(\psi)$ and minimize w.r.t. $R$, obtaining
an upper bound on $\Ffrak(\psi_1)$.
 Rescaling to units in which the ground state energy $E=-1$ yields r.h.s.\refeq{normOFuSQRbounds}.
\hfill $\square$

 We end this section by remarking that our lower and upper bounds \refeq{normOFuSQRbounds} are 
slightly stronger than those one can obtain from the energy bounds in \cite{T} by scaling.

\section{Numerical illustration of the Schr\"odinger--Newton ground state}\vspace{-10pt}

\begin{figure}[H]
 \includegraphics[width = 10truecm,height=12truecm]{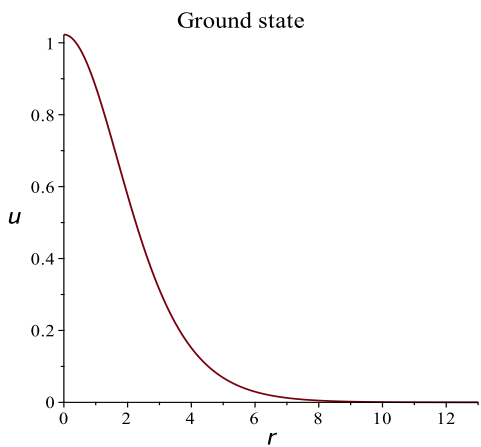}
\hspace{-.7truecm}
\vspace{-6truecm}
\caption{Shown is a numerical approximation to the gound state solution $u(r)$ of \refeq{CPSN}.}   
\label{SN}
\end{figure}

\begin{figure}[H]
  \includegraphics[width = 10truecm,height=12truecm]{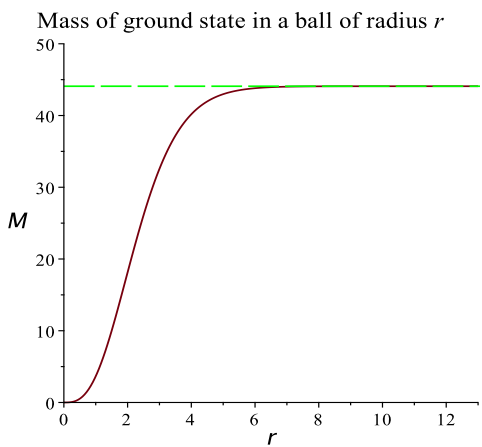} 
 \vspace{-5.7truecm}
\caption{Shown is $M(r):=4\pi\int_0^r |u(s)|^2s^2ds$ for the gound state solution $u(r)$ of \refeq{CPSN}, together with
the horizontal asymptote at $14.03\pi\approx \|u\|_2^2 =\lim_{r\to\infty}M(r)$.}   \vspace{-.5truecm}
\label{SNmassVSr}
\end{figure}

We have numerically computed an approximation to 
the Schr\"odinger--Newton ground state by dividing \refeq{CPSN} by $u(r)$, then
applying $-\Delta$, which by the radial symmetry yields a fourth-order ordinary differential equation for $u(r)$,
for which initial data $u(0)$, $u^\prime(0)$, $u^{\prime\prime}(0)$, and $u^{\prime\prime\prime}(0)$ need to be supplied. 
 From \refeq{CPSN} one inherits $u^\prime(0)=0$, and 
well-posedness of the initial value problem for the fourth-order ODE requires $u^{\prime\prime\prime}(0)=0$ also. 
 This leaves $u(0)>0$ and $u^{\prime\prime}(0)$ to be determined such that $u(r)>0$ for all $r\geq 0$, 
with $\ln u(r) \sim -r$ as $r\to\infty$, and $u\in L^2(\Rset,r^2dr)$. 
 This single fourth-order ODE formulation is equivalent to the system of two second-order ODEs discussed in
\cite{TM}, \cite{MPT}, \cite{Hthesis}, and \cite{GW}.
   Fig.~\ref{SN} shows the ground state $u(r)$ and Fig.~\ref{SNmassVSr} its mass function $M(r):=4\pi\int_0^r |u(s)|^2s^2ds$.
 Note that $\|u\|_2^2 =\lim_{r\to\infty}M(r)$.
\vspace{-1truecm}
 \vspace{-1.5truecm}

\begin{figure}[H]
\hspace{-.5truecm}
  \includegraphics[width = 10.5truecm,height=11truecm]{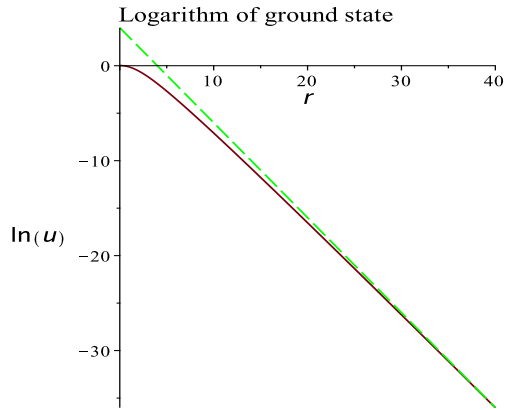} \vspace{-5truecm}
\caption{Shown is the natural logarithm of the gound state solution $u(r)$ of \refeq{CPSN}, together with
a straight line of slope $-1$. }  \vspace{-1.7truecm}
\label{SNlog}
\end{figure}

\begin{figure}[H]
\hspace{-.5truecm}
  \includegraphics[width = 10.5truecm,height=12truecm]{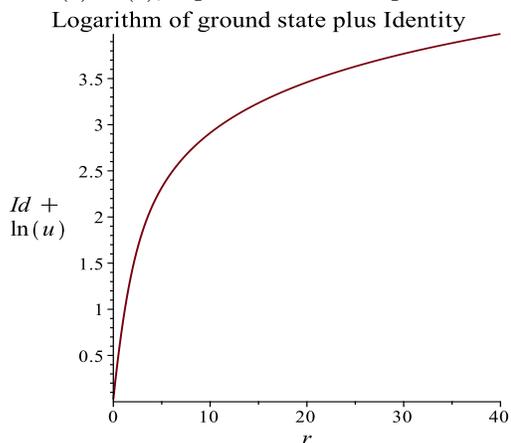} \vspace{-6truecm}
\caption{Shown is $r+ \ln u(r)$ versus $r$ for the gound state solution $u(r)$ of \refeq{CPSN}.}  
\label{SNlogOfuPLUSid} \vspace{-.5truecm}
\end{figure}

Fig.~\ref{SNlog}, where we display the natural logarithm of $u(r)$ together with a straight line of
slope $-1$, appears to suggest a purely exponential decay of the ground state $u(r)$, but appearances are
misleading, as visualized in Fig.~\ref{SNlogOfuPLUSid}.

$\phantom{nix}$\vspace{-1.6truecm}

\begin{figure}[H]
\hspace{-.5truecm}
  \includegraphics[width = 11truecm,height=13truecm]{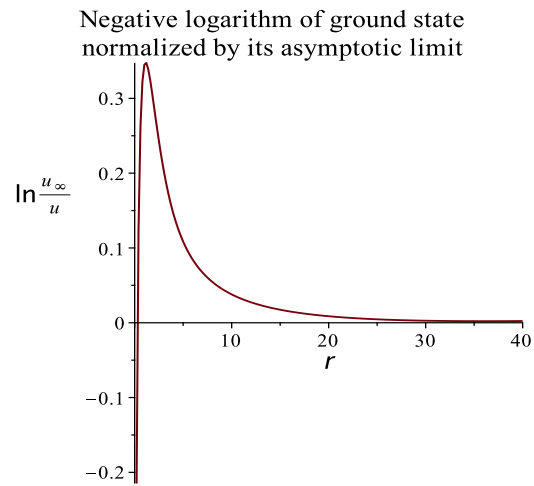} \vspace{-6truecm}
\caption{Shown is the negative natural logarithm of the ratio of the gound state solution $u(r)$ of \refeq{CPSN} 
over its asymptotic limit $u_\infty(r) = A r^\beta\exp(-r)$ with $\beta \approx 0.754$ and $A\approx 3.37$.}  \vspace{-.5truecm}
\label{SNlnRATIOuOVERuVSr}
\end{figure}

Fig.~\ref{SNlogOfuPLUSid} reveals that the map $r\mapsto r + \ln u(r)$ is not asymptotic, for large $r$, to
a constant function, which it would be if $u(r)\sim A \exp(-r)$. 
 Instead, this map seems to behave $\propto \ln r$ for large $r$, which is confirmed in Fig.~\ref{SNlnRATIOuOVERuVSr}.
 \vspace{-1.5truecm}

 Our numerical computations were carried out with MAPLE's 
Cash--Karp fourth-fifth order Runge--Kutta method with degree four interpolant (ck45), which proved
more suitable than MAPLE's default Runge--Kutta--Fehlberg routine rkf45. 
 To overcome the enormous variations over the range of $u(r)$ we solved the ODE for $\ln u(r)$ 
and asked for 70 digits precision during the computation.
 The interval halving iterations to determine the correct initial data $u(0)$ and $u''(0)$ to 
yield $\ln u(r) \sim -r$ were terminated after a precision of three significant digits had been achieved, though.
 
 As a test for our results we rescaled the ground state energy $E=E_0$ for \refeq{CPSNnat}, computed 
numerically in \cite{GW} to be $E_0=-0.325(74)$, into our units in which the ground state energy $E=-1$.
 This yields $\|u\|_2^2 = 8\pi/\sqrt{|E_0|}\approx 44.09$, in good agreement with our result $14.03\pi\approx 44.08$. 
\vspace{-20pt}

\section{External $\sim - K/r$ potentials}\vspace{-10pt}

 In its bosonic star interpretation, the Schr\"odinger--Newton equation \refeq{CPSN} 
captures the quantum-mechanical ground state of a single species of $N$ spin-zero bosons that interact only with 
Newtonian gravity among themselves, in the Hartree limit $N\to\infty$; cf. \cite{Lewin}, and see \cite{LiebYauCMP} 
for a version with special-relativistic kinetic energy.
 Since these hypothetical building blocks of the mysterious dark matter will interact gravitationally
with other matter, it is of interest to also consider the Schr\"odinger--Newton equation for bosons which
are exposed to an external gravitational potential, which for simplicity we will assume to be radially 
symmetric  and asymptotically $\sim - K/r$ for some $K>0$.
 Incidentally, equipped with such an additional external $\sim - K/r$ potential the so-generalized
Schr\"odinger--Newton equation also can have bound states if the sign of the self-consistent interaction
term (the cubic nonlinearity term at r.h.s.\refeq{CPSN}) is changed from ``$+$'' to ``$-$.'' 
 Therefore, in this and the next section we consider the equation
\begin{eqnarray}\hspace{-20pt}
-  \Delta u(\sV) + u(\sV) =
 \displaystyle{\frac{1}{4\pi}\int_{\Rset^3} \!\!\!\frac{\mathfrak{z}\|u\|_2^2\mu(\sV')\pm |u|^2(\sV')}{\abs{\sV-\sV'}} d^3\!s'} u(\sV)
 \label{CPSNH} 
\end{eqnarray}
with $\mathfrak{z}>0$ if the ``$+$'' sign is chosen, and in this section with $\mathfrak{z}> 1$ if the ``$-$'' sign is;
in the next section we also allow $\mathfrak{z}\geq 0.825$.
 Here, $\mu(\sV)\geq 0$ is a given probability measure, in the following assumed to be radially symmetric.
 Eq.\refeq{CPSNH} with the ``$-$'' sign and $\mathfrak{z}>1$ captures the ground state of a
positive ion with hypothetical bosonic ``electrons'' in the Hartree limit $N\to\infty$,  
\cite{BenguriaLieb}, \cite{LiebSeiringer}; cf. also \cite{BEGMY}, \cite{FKS}, \cite{P}, \cite{KieJMP}, \cite{LNR}, \cite{Rou},
and references therein. 
 In that case $\mu(\sV)$ represents the (normalized) charge distribution of the nucleus, which usually
is modelled by a Dirac $\delta$ measure but may also be considered to be bounded and radially symmetric decreasing.
 In the gravitational interpretation (i.e. ``$+$'' sign in \refeq{CPSNH}),  $\mu(\sV)$ is the normalized mass distribution
of the other matter. 

 Eq.\refeq{CPSNH}, with either sign, is the formal Euler--Lagrange equation for the minimization of the functional
\vspace{-10pt}
\begin{widetext}
\begin{eqnarray} \label{eq:uFctlEXT}\hspace{-10pt}
\Ffrak^\mp(u) := 
\int_{\Rset^3} \left(\abs{\nabla u}^2
+
 |u|^2\right)\!(\sV)\, d^3\!s 
\mp \frac{1}{8\pi} \int_{\Rset^3} \int_{\Rset^3} 
\frac{|u|^2(\sV)|u|^2(\sV')}{\abs{\sV-\sV'}}d^3\!s d^3\!s'
- \mathfrak{z}\frac{\|u\|_2^2}{4\pi} \int_{\Rset^3} \int_{\Rset^3} 
\frac{|u|^2(\sV)\mu(\sV')}{\abs{\sV-\sV'}}d^3\!s'd^3\!s 
\end{eqnarray}
\end{widetext}
$\phantom{nix}$\vspace{-25pt}

\noindent
over the set of $u\in H^1(\Rset^3)$ under the Nehari condition that $u\neq 0$ satisfies
${4\pi} \int_{\Rset^3}\big( \abs{\nabla u}^2(\sV) + |u|^2(\sV) \big) d^3s 
= 
\int_{\Rset^3} \int_{\Rset^3} {\abs{\sV-\sV'}}^{-1}{|u|^2(\sV)\big[\mathfrak{z}\|u\|_2^2\mu(\sV')\pm|u|^2(\sV')\big]}d^3s d^3s'$;
no normalization condition on $\int |u|^2 (\sV)d^3s$ is imposed.
 Under our assumptions that  $\mathfrak{z} >0$ if the upper sign is chosen, and $\mathfrak{z}> 1$ if the lower
one is, we may expect that a unique positive minimizer of $\Ffrak^\mp(u)$ exists and satisfies \refeq{CPSNH},
and in the following we assume this; cf.~\cite{KK}.

 The asymptotic behavior of such a positive solution of
\refeq{CPSNH} for rapidly decreasing $\mu(r)$ (bounded by a Schwartz function) can be
obtained by adapting the strategy of  V. Moroz and J. van Schaftingen \cite{MvS} to  \refeq{CPSNH}. 
 We find 
\begin{equation}\label{asympCPSNHext}
\lim_{r\to\infty} r^{1 - (\mathfrak{z}\pm 1)\|u\|_2^2/8\pi} e^r u(r) \in (0,\infty);\vspace{-10pt}
\end{equation}
thus, $u(r) = A e^{-r}/ r^{1 - (\mathfrak{z}\pm 1)\|u\|_2^2/8\pi}+$l.o.t. for an $A>0$. 

\vspace{-17pt}

\section{Concluding remarks}\vspace{-10pt}
 In section IV we have only considered the Schr\"odinger--Newton and Hartree equations in settings which allow
a straightforward adaptation of the strategy of \cite{MvS}. 
 In the ionic Hartree problem with $\mu(\sV) = \delta(\sV)$, it is known that a minimizer 
exists also for $\mathfrak{z}_* \leq \mathfrak{z} \leq 1$ with $\mathfrak{z}_*\approx 0.825$; 
cf.\cite{BenguriaLieb}, \cite{LiebSeiringer}, \cite{Baum}.
 Although we did not investigate this more difficult neutral atom and
negative ion regime to the point that we can make a definite
statement, our preliminary inquiry leads us to the surmise that
\refeq{asympCPSNHext} continues to hold for neutral bosonic atoms and negative ions, so that in this regime one seems to have
$u(r)\sim A e^{-r}\!/r^\gamma$ with $\!\gamma \geq 1$.
  Interestingly, \refeq{asympCPSNHext} thus suggests that
the asymptotic form of $u(r)$ proposed by K. P. Tod and I. Moroz in \cite{TM} for
\refeq{CPSN} is instead true for the Hartree ground state of a neutral 
bosonic atom in the $N\to\infty$ limit, when $\mathfrak{z}=1$. 
 In all the cases discussed in section IV the asymptotic form $u(r)\sim A e^{-r}/r^\gamma$ holds with $\gamma<1$.
 
 Given that the question of the asymptotic decay of the Schr\"odinger--Newton ground state
has received several conflicting answers until rigorous results showed the way, it is desirable
to rigorously vindicate also the findings and the surmise reported in sections IV and V.
 \vspace{-17pt}

\section*{Acknowledgments} \vspace{-12pt}
 I thank Parker Hund and Eric Ling for discussions. \vspace{-5pt}
\medskip

\textbf{Note added:} Corrections (\textcolor{red}{red}) and minor editorial additions (\textcolor{blue}{blue})
on March 06, 2021. A few now obsolete statements in the published version (Phys. Lett. A \textbf{395}, art. 127209 (2021))
have been omitted.

\newpage

\section*{\hspace{-9truecm}References} \vspace{-20pt}  

\bibliographystyle{apsrev}

\medskip\medskip

miki@math.rutgers.edu

\end{document}